\documentstyle[12pt]{article}
  \newcommand{\const}{\rm const}

 \newcommand{\sign}{\rm sign}
 \newcommand{\mes}{\rm mes}

\begin{document}

   \begin{center}

  \   {\bf  Exact asymptotic for tail of distribution of self-normalized }\\

\vspace{4mm}

   \  {\bf sums of random variables under classical norming. }\\

\vspace{4mm}

  \  {\bf Ostrovsky E., Sirota L.}\\

\vspace{4mm}

 Israel,  Bar-Ilan University, department of Mathematic and Statistics, 59200, \\

\vspace{4mm}

E-mails: \\
 eugostrovsky@list.ru,  \ sirota3@bezeqint.net \\

\vspace{5mm}

  {\bf Abstract} \\

\vspace{4mm}

 \end{center}

\ We derive in this article the  asymptotic behavior  as well as non-asymptotical estimates
of tail of distribution  for self-normalized sums  of  random variables (r.v.) under natural classical norming. \par

 \ We investigate also the case of non-standard random norming function and the tail asymptotic for the maximum distribution
for self-normalized statistics. \par

 \ We do not suppose the independence or identical distributionness of considered random variables, but we assume the
existence and sufficient smoothness of its density.\par

 \ We  show also the exactness  of our conditions imposed on the considered random variables  by means of  building of an
appropriate examples (counterexamples). \par

\vspace{7mm}

 \par

\vspace{4mm}


\vspace{5mm}

\section{ Definitions.  Notations. Previous results.  Statement of problem.}

 \vspace{3mm}

 \ Let  $ \ \{ \xi(i) \}, \ i =1, 2, 3,\ldots,n; \ n \ge 2, \  $ be a collection of random variables or equally random vector (r.v.)

$$
  \xi = \vec{\xi} =  \{  \xi(1), \xi(2), \ldots, \xi(n)  \},
$$
 not necessary to be independent, centered or identically distributed,  defined on
certain probability space,  $  \ \forall i \ \Rightarrow {\bf P} ( \ \xi(i) = 0)  = 0, \ $
having a (sufficiently smooth) density of distribution $  \  f_{\vec{\xi}}(\vec{x}) =  f(x) = f(\vec{x}), \ x = \vec{x} \in  R^n.   \ $ \par
 \ Let us introduce the following self-normalized sequence of sums of r.v. under the classical norming

$$
 T = T(n) =  \frac{\sum_i \xi(i)}{ \sqrt{\sum_i \xi^2(i) }}, \eqno(1.1)
$$
here and in what follows

$$
\sum = \sum_i = \sum_{i=1}^n, \ \sum_j = \sum_{j=2}^n, \ \prod_j = \prod_{j=2}^n,
$$
and define the correspondent tail probabilities

$$
Q_n = Q_n(B) := {\bf P}(T(n) > B), \ B = \const > 0;
$$

$$
Q(B) := \sup_n Q_n(B)  = \sup_n  {\bf P}(T(n) > B), \ B = \const > 0.
$$

 \ B.Y.Jing, H.Y.Liang and W.Zhou obtained in an article [5] the following uniform estimate  of sub-gaussian type for i, i.d.  symmetrical
non-degenerate random variables  $  \ \{ \xi(i) \} \  $

$$
Q(B)  \le \exp  \left(  - B^2/2 \right).
$$

 \  Note first of all that if $  \ n = 1, \ $ then the r.v. $ \ T(1)  = \sign(\xi(1)) $  has a Rademacher's distribution. This case is trivial for us
and may be excluded.  \par

 \ Further, it follows from the classical  H\''older's inequality that $ \ T(n) \le \sqrt{n}, \ n = 2,3, \ldots; $ therefore

$$
 \forall B \ge \sqrt{n} \Rightarrow  Q_n(B) = 0.
$$
 \ Thus, it is reasonable to suppose $  \ B = B( \epsilon )  = \sqrt{n} -\epsilon, \ \epsilon \in (0,1); $ and  it is interest by our
opinion to investigate the asymptotical behavior  as well as the non-asymptotical estimates  for  the following  tail function

$$
q(\epsilon) = q_n(\epsilon) = Q_n(\sqrt{n} - \epsilon) \eqno(1.2)
$$
as $  \ \epsilon \to 0+, \ \epsilon \in (0,1); $  the value $  \ n \ $ will be presumed to be fix and greatest or equal than 2. \par

 \ The case of the left tail of distribution $  \ {\bf P} (T(n) < - \sqrt{n} + \epsilon), \  \epsilon \in (0,1), \ $
as well as the probability $  \ {\bf P} (|T(n)| > \sqrt{n} - \epsilon), \  \epsilon \in (0,1), \ $
may be investigated quite  analogously. \par

\vspace{4mm}

  \ {\bf Our purpose in this short preprint is just obtaining an  asymptotical expression of these probabilities,
as well as obtaining the non-asymptotical bilateral estimates for ones. } \par

 \ {\bf   We consider also the case of non-standard random norming function.  } \par

\vspace{4mm}

 \  The problem of tail investigation for self-normalized  random sums
with at the same or another self norming sequence was considered in many works,
see e.g. [1]-[11].  Note that in these works
was considered as a rule only asymptotical approach, or uniform estimates,  i.e. when $ \ n \to \infty; \ $  for
instance, was investigated the classical Central Limit Theorem (CLT), Law of
Iterated Logarithm (LIL) and Large Deviations (LD) for these variables.
Several interest applications of these tail functions, in particular,
in the non-parametrical statistics are described in [1], [2], [5], [7]-[8], [11] etc. \par

\vspace{4mm}

\section{ Main result. }

 \vspace{4mm}

 \ We need to introduce now some needed notions and notations.  Introduce  for any $ \ n \ $ dimensional vector
$  \ x = \vec{x}= \{x(1), x(2), \ldots, x(n)  \}  \ $ its $ (n-1) \ - $ dimensional sub-vector

$$
y = \vec{y} = \vec{y}(\vec{x}) :=  \{x(2), x(3), \ldots, x(n)\}. \eqno(2.1)
$$

 \ Define also for arbitrary $ \ (  n-1 )  - \ $ dimensional positive  vector  $ \ v = \vec{v} = \{ v(2), v(3), \ldots, v(n)  \}   \ $  the function

$$
g(v) = g_n(v) = \frac{1 + \sum_j v(j)}{\sqrt{ 1 + \sum_j v^2(j)  }}  \eqno(2.2)
$$
and introduce the correspondent its {\it anti-Hessian} matrix for this function at the extremal point
$  \  \vec{v_0} = \vec{1} =  (1,1, \ldots,1),  \  \dim \vec{v_0} = (n-1), \ $ containing the following entries: $ \ A =  A(n-1) = \{  a(j,k) \},  $

$$
 a(j,k)  := - \left\{  \frac{\partial^2 g(v) }{\partial v(j) \ \partial v(k)}  \right\}/ \vec{v} = \vec{1}, \ j,k = 2,3,\ldots,n; \eqno(2.3)
$$
  and we find by the direct computations

$$
a(j,j)  = n^{-1/2} - n^{-3/2} - \eqno(2.3a)
$$
the diagonal members,

$$
 \ a(j,k) = - n^{-3/2}, \ k \ne j -  \eqno(2.3b)
$$
off diagonal entries. \par

\vspace{4mm}

 \ {\bf Lemma 2.1.}  Let $  \  L_m = L_m(x), \  m = 1,2,3, \ldots \ $ be a square matrix of a size $ \ m \times m $ with entries

$$
l(j,j) = x, \ x \in R; \  l(j,k) = 1; \ j,k = 1,2,\ldots,m; \ j \ne k.
$$
 \ Then

$$
\det L_m = (x-1)^{m-1} \cdot (x - m + 1).
$$

\vspace{4mm}

 \ {\bf  Corollary 2.1.}  Let $  \  L_m = L_m(a,b), \  m = 1,2,3, \ldots \ $ be a square matrix of a size $ \ m \times m $ with
entries

$$
l(j,j) = a,  \  l(j,k) = b, \ a,b \in R; \ j,k = 1,2,\ldots,m; \ j \ne k.
$$
 \ Then

$$
\det L_m(a,b) =  (a -b )^{m-1} \cdot (a - (m - 1)b).
$$

\vspace{4mm}

 \  It is no hard to compute by virtue of Corollary 2.1 the determinant of the introduced before matrix $  \ A, \ $ which will be
used further:

$$
\det(A) =  n^{-(n-2)/2}  \cdot \left(2n^{-1/2} - 3 n^{-3/2} \right). \eqno(2.3c).
$$

 \ Note that this matrix  $ \ A \  $ is symmetric and positive definite. \par

 \ Further, define a following function as an integral

$$
h(\vec{v})  :=  \int_{-\infty}^{\infty}  \ f(z,  \ z \ \vec{v} )  \ dz, \eqno(2.4)
$$
so that

$$
h(\vec{1}) =  \int_{-\infty}^{\infty} f_{\vec{\xi}} (z, z, \ldots, z) \ dz. \eqno(2.4a)
$$

\vspace{4mm}

{\bf Theorem 2.1.} Suppose that the function $   \ h(\vec{v}) \  $ there exists,   $   \ h(\vec{1}) > 0, \  $
 and is continuous at the point $  \ \vec{v_0} = \vec{1}. \  $ Then for  (positive finite) constant $  K = K(n): $

$$
K(n) :=  2^{-(n-1)/2} (\det A(n-1))^{-1/2} \ \frac{\pi^{ (n-1)/2 }}{\Gamma( (n+1)/2 )} =
$$

$$
2^{-(n-1)/2} \ n^{(n - 2)/4} \  \left(2n^{-1/2} - 3 n^{-3/2} \right)^{-1/2} \  \frac{\pi^{ (n-1)/2 }}{\Gamma( (n+1)/2 )}.
 \eqno(2.5)
$$
we have

$$
q_n(\epsilon) \sim  K(n) \ h(\vec{1}) \  \epsilon^{(n-1)/2}.  \eqno(2.6)
$$

\vspace{4mm}

 \ {\bf Proof.} Note first of all that the point $  \vec{v_0} = \vec{1} $ is an unique point of maximum of the smooth function
$ \ v \to g(v); $ and this maximum is equal to $ \sqrt{n}. $ \par

\ Further, we have as $  \ \epsilon \to 0+ $

$$
q_n(\epsilon) = {\bf P} \left(  \frac{\sum \xi(i)}{ \sqrt{ \sum \xi^2(i)  }  }  > \sqrt{n}  - \epsilon  \right) =
$$

$$
\int \int \ldots \int_{  \sum x(i) /\sqrt{ x^2(i)} > \sqrt{n}  - \epsilon }   f(x(1),  x(2), \ldots, x(n)) \ dx(1) \ dx(2), \ldots \ dx(n)  =
$$

$$
\int \int \ldots \int_{  \left[x(1) + \sum y(j) \right]/\sqrt{ x^2(1) + \sum_j y^2(j)} > \sqrt{n} - \epsilon  } \  \cdot \
f(x(1), \vec{y}) \ dx(1) \ dy =
$$

$$
 \int_0^{\infty} dx(1) \int_{ [1 + \sum_j v(j)] /\sqrt{ 1 + \sum_j v^2(j)   }  > \sqrt{n}  - \epsilon } \ \prod_j v(j) \ \cdot  f(x(1), x(1)\vec{v} ) \ d \vec{v}   =
$$

$$
 \int_0^{\infty} dx(1) \int_{g(\vec{v})  > \sqrt{n} - \epsilon } \ \prod_j v(j) \cdot \  f(x(1), x(1) \ \vec{v} ) \ d \vec{v}   =
$$

$$
 \int_{g(\vec{v})  > \sqrt{n} - \epsilon }  \ \prod_j v(j) \cdot \ h(\vec{v}) \ d \vec{v} =
 \int_{g(\vec{v})  > \max g(\vec{v}) - \epsilon }  \ \prod_j v(j) \cdot \ h(\vec{v}) \ d \vec{v}  =
$$

$$
 \int_{g(\vec{v})  > g(\vec{1}) - \epsilon }  \ \prod_j v(j) \cdot \ h(\vec{v}) \ d \vec{v}.  \eqno(2.7)
$$
  \ The last integral is localized in some sufficiently small neighborhood of the point of maximum $ \vec{v} = \vec{v_0} = \vec{1}.  $  In detail,
 notice that as $  \ \epsilon \to 0+ \ $ the set $  \  \{  v: \  g(\vec{v})  > g(\vec{1}) - \epsilon   \}  \  $ is asymptotical equivalent to the ellipsoidal
set

$$
\  \{  v: \  (A(v-1), (v-1)) < \epsilon \},
$$
therefore

$$
q_n(\epsilon) \sim \int_{ \{  v: \  (A(v-1), (v-1)) < \epsilon \}   }  \ \prod_j v(j) \cdot \ h(\vec{v}) \ d \vec{v},
$$
which is in turn asymptotical equivalent to the following integral

$$
q_n(\epsilon) \sim   h(\vec{1}) \cdot   \int_{ \{  v: \  (A(v-1), (v-1)) < \epsilon \}   }  \ d \vec{v} =
\mes  \{  v: \  (A(v-1), (v-1)) < \epsilon \},
$$
 and we find after simple calculations

$$
q_n(\epsilon) \sim  h(\vec{1}) \ K(n) \cdot \epsilon^{(n-1)/2}, \eqno(2.8)
$$
Q.E.D. \par

\vspace{4mm}

 \ {\bf Remark 2.1.} If the function $  \ v \to h(v) \ $ is  not continuous but only  integrable in some sufficiently small neighborhood
of the point $  \ \vec{1}, \ $  or perhaps

$$
\lim_{ ||\vec{v} - \vec{1}|| \to 0   }   \prod_j v(j)  \cdot   h(\vec{v}) = 0,
$$
then

$$
q_n(\epsilon) \sim    \int_{ \{  v: \  (A(v-1), (v-1)) < \epsilon \}   } \prod_j v(j) \cdot  h(\vec{v}) \ d \vec{v},
$$
if of course the last integral is finite and non-zero. \par

  \ Assume for instance that for $  \vec{v} \to \vec{1}  $

$$
\prod_j v(j) \cdot  h(\vec{v})  \sim  \left[(A(v-1), (v-1))  \right]^{\gamma/2}, \ \gamma = \const > 1 - n;
$$
then  as $  \ \epsilon \to 0+ $

$$
q_n(\epsilon) \sim 2^{ -(n-3)/2  } \cdot
(\det A)^{-1/2} \cdot \frac{\pi^{(n-1)/2}}{\Gamma((n-1)/2)} \cdot \frac{\epsilon^{ ( n + \gamma - 1  )/2  }}{n + \gamma - 1}.
$$

\vspace{6mm}

 \  Let us return to the  promised above case of the left tail of distribution $  \ {\bf P} (T(n) < - \sqrt{n} + \epsilon), \
\epsilon \in (0,1), \ $  as well as the case of the probability $  \ {\bf P} (|T(n)| > \sqrt{n} - \epsilon), \  \epsilon \in (0,1). \ $ \par

\vspace{4mm}

{\bf  Corollary 2.1.} Suppose that the function $   \ h(\vec{v}) \  $ there exists,   $   \ h( - \vec{1}) > 0, \  $
 and is continuous at the point $  \ \vec{v_-} = - \vec{1}. \  $ Then for  at the same positive finite constant $  K = K(n) $
we have

$$
  \ {\bf P} (T(n) < - \sqrt{n} + \epsilon)  \sim  K(n) \ h( - \vec{1}) \  \epsilon^{(n-1)/2}.  \eqno(2.9)
$$

\vspace{4mm}

{\bf  Corollary 2.2.} Suppose that the function $   \ h(\vec{v}) \  $ there exists,   $  \ h(\vec{1}) + \ h( - \vec{1}) > 0, \  $
 and is continuous at  both the the points $ \  \vec{v_0} =\vec{1} \ $ and $  \ \vec{v_-} = - \vec{1}. \  $ Then for  at the same
positive finite constant $  K = K(n) $ we have  as $ \  \epsilon \to 0+ \  $

$$
  \ {\bf P} (|T(n)| > \sqrt{n} - \epsilon)  \sim  K(n) \  [ h(\vec{1}) + h( - \vec{1}) ] \  \epsilon^{(n-1)/2}.  \eqno(2.10)
$$

\vspace{4mm}

\section{ Some generalizations: non-standard norming random function.}

\vspace{4mm}

 \ Let $ \  \beta = \const > 1   \  $ and the sequence of r.v. $  \  \{\xi(i) \}, \ i = 1,2,\ldots, n; \  n \ge 2 \ $ is as before.
The following statistics was introduced (with applications) at first perhaps by Xiequan Fan [2]:

$$
T_{\beta}(n) \stackrel{def}{=}  \frac{\sum \xi(i)} {  [ \sum |\xi(i)|^{\beta}   ]^{1/\beta} }. \eqno(3.1)
$$

\ Xiequan Fan derived  in particular in [2] the following generalization of result belonging to
 B.Y.Jing, H.Y.Liang and W.Zhou \  [5] \   of sub-gaussian type for i, i.d. symmetrical
non-degenerate r.v. $  \ \{ \xi(i) \} \  $

$$
 {\bf P} \left( T_{\beta}(n)   > B \right) \le \exp \left( - 0.5 \ B^2 \ n^{ 2/\beta - 1  }   \right),  \ \beta \in (1,2].
$$

 \ Since the theoretical  attainable  maximum of these statistics is  following:

$$
\sup_{ \{\xi(i) \}} T_{\beta}(n) = n^{1 - 1/\beta},
$$
it is reasonable to investigate the next tail probability

$$
r_{\beta, n} (\epsilon) = r(\epsilon)  := {\bf P} \left( T_{\beta}(n) >   n^{1 - 1/\beta} - \epsilon \right), \ \epsilon \to 0+, \epsilon \in (0,1). \eqno(3.2)
$$

 \ Define the following modification of the $  \ g \ -  $ function:

$$
g_{\beta}(\vec{v}) := \left( 1 + \sum_j v_j \right) \cdot  \left( 1 + \sum_j v_j^{\beta} \right)^{-1/\beta}, \  \vec{v} \in R^{n-1}_+ \eqno(3.3)
$$
which attained its maximal value  as before at the point  $ \  v = \vec{v} = \vec{1}  \  $ and herewith

$$
\max_v g_{\beta}(\vec{v}) = g_{\beta}(\vec{1}) = n^{1 - 1/\beta}, \eqno(3.4a)
$$

$$
-\frac{\partial g_{\beta}^2}{\partial v_k \ \partial v_l}(\vec{1}) = - (\beta - 1) n^{ -1 - 1/\beta  }, \ k \ne l; \eqno(3.4b)
$$

$$
-\frac{\partial g_{\beta}^2}{\partial v_k^2}(\vec{1}) = (\beta - 1) \left[ n^{ - 1/\beta  }  - n^{ - 1 - 1/\beta  } \right]. \eqno(3.4c)
$$

 \ The correspondent {\it anti-Hessian} matrix  $  \ A_{\beta} = A_{\beta}(n -1) $
 at the same  extremal point $  \  \vec{v_0} = \vec{1} =  (1,1, \ldots,1),  \ \dim \vec{v_0} = \ (n-1) \ $ contains the following entries:
$ \ A_{\beta} =  A_{\beta}(n-1) = \{  a_{\beta}(j,k) \},  $  where

$$
 a_{\beta}(j,k)  := - \left\{  \frac{\partial^2 g_{\beta}(v) }{\partial v(j) \ \partial v(k)}  \right\}/ \vec{v} = \vec{1}, \ j,k = 2,3,\ldots,n; \eqno(3.5)
$$
  and we find  by direct computations

$$
a_{\beta}(j,j)  = (\beta - 1) \ \left[ \ n^{-1/\beta} - n^{- 1 - 1/\beta} \ \right] \ - \eqno(3.6a)
$$
diagonal members,

$$
 \ a_{\beta}(j,k) =  -(\beta - 1)  \ n^{- 1 - 1/\beta}, \ k \ne j   -  \eqno(3.6b)
$$
off diagonal entries. \par

 \ The correspondent determinant has a form

$$
 \det A_{\beta}(n-1) = (\beta - 1)^{n-1} \cdot n^{ - (n-2)/\beta }  \cdot \left[  2 n^{-1/\beta} - 3 n^{  - 1 - 1/\beta  }   \right]. \eqno(3.7)
$$

 \ Of course, the last expressions (3.5)-(3.7) coincides with ones when $ \ \beta = 2 \ $ in the second section.  \par
 \ We deduce  similar to the second section

\vspace{4mm}

{\bf Theorem 3.1.} Suppose  as before that the function $   \ h(\vec{v}) \  $ there exists,   $   \ h(\vec{1}) > 0, \  $
 and is continuous at the point $  \ \vec{v_0} = \vec{1}. \  $ Then for  (positive finite) constant $  K_{\beta} = K_{\beta}(n): $

\vspace{3mm}

$$
K_{\beta}(n) =  2^{-(n-1)/2} \ (\det A_{\beta}(n-1))^{-1/2} \ \frac{\pi^{ (n-1)/2 }}{\Gamma( (n+1)/2 )} =
$$

$$
2^{-(n-1)/2} \  n^{ ( n -2)/(2 \beta) }  \cdot \left[  2 n^{-1/\beta} - 3 n^{  - 1 - 1/\beta  }   \right]^{-1/2} \cdot  \frac{\pi^{ (n-1)/2 }}{\Gamma( (n+1)/2 )}
\eqno(3.8)
$$

we have as  $ \ \epsilon \to 0+  $

$$
r_{\beta, n} (\epsilon)   \sim  K_{\beta}(n) \ h(\vec{1}) \  \epsilon^{(n-1)/2}.  \eqno(3.9)
$$

 \vspace{4mm}

\ {\bf Remark 3.1.}  It follows immediately from Stirling's  formula  that as $ \ n \to \infty   \ $

$$
\log_n K_{\beta}(n) \sim  n \cdot \frac{1 - \beta}{2 \beta}. \eqno(3.10)
$$

 \ Thus, the sequence $ \  K_{\beta}(n)  \  $ tends as $ \ n \to \infty \ $ very rapidly to zero. Recall that we consider the case when
$  \ \beta > 1. $ \par

\vspace{4mm}

  \section{ Tail of maximum distribution estimates.}

\vspace{4mm}

 \  Define following Xiequan Fan [2]  the  tails of maximum distributions

$$
R_n(\epsilon) := {\bf P} \left(  \max_{k=2,3,\ldots,n} \frac{ S(k) }{Z(n)} > \sqrt{n} - \epsilon  \right), \eqno(4.1)
$$
where

$$
S(k) = \sum_{l=1}^k \xi(i), \ Z(n) = \sqrt{ \sum_j \xi^2(j)   },
$$
and

$$
 \epsilon \in \left( \ 0, [  2 \sqrt{n-1}  ]^{-1} \right). \  \eqno(4.2)
$$
 \ We aim to investigate as before the asymptotic behavior as  $ \ \epsilon \to 0+ \  $ of this tail probability. \par

\vspace{4mm}

 \ {\bf Theorem 4.1.}  We propose under at the same conditions as in theorem 2.1 that as $  \ \epsilon \to 0+ $ and
$  \ \epsilon \in \left( \ 0, [  2 \sqrt{n-1}  ]^{-1} \right). \ $

$$
R_n(\epsilon) \sim Q_n(\epsilon).  \eqno(4.3)
$$

\vspace{4mm}

 \ {\bf Proof.}  The lower bound is trivial:

$$
R_n(\epsilon) := {\bf P} \left(  \max_{k=2,3,\ldots,n} \frac{ S(k) }{   Z(n)} > \sqrt{n} - \epsilon  \right) \ge
{\bf P} \left( \frac{ S(n) }{   Z(n)} > \sqrt{n} - \epsilon      \right) = Q_n(\epsilon).
$$
 \ It remains to ground the inverse inequality. One can suppose without loss of generality $  \ \xi(i) > 0. \  $ We have:

$$
R_n(\epsilon) \le \sum_{k=2}^n R_{n,k}(\epsilon), \eqno(4.4)
$$
where

$$
R_{n,k}(\epsilon) = {\bf P} \left(  \frac{S(k)}{ Z(n) }   > \sqrt{n} - \epsilon \right).
$$

 \ Let now $ \  2 \le k \le n - 2;  \ $ then

$$
 \frac{S(k)}{ Z(n) } \le  \frac{S(k)}{ Z(k) } \le \sqrt{k} < \sqrt{n} - \epsilon,
$$
\ therefore

$$
 R_{n,k}(\epsilon) = 0,
$$
if $ \ \epsilon \  $ satisfoes the restriction (4.2). Thus,

$$
R_n(\epsilon) = R_{n,n}(\epsilon) = Q_n(\epsilon),
$$
Q.E.D. \par

\vspace{4mm}

 \ Note in addition that we have proved in fact that

$$
\overline{R}_n(\epsilon) := {\bf P} \left(  \max_{k=2,3,\ldots,n} \frac{ S(k) }{Z(k)} > \sqrt{n} - \epsilon  \right)
\sim Q_n(\epsilon), \ \epsilon \to 0+. \eqno(4.5)
$$

\vspace{4mm}

 \section{ Non-asymptotical estimates.  }

\vspace{4mm}

  \ Introduce the following important functional

$$
\lambda = \lambda(g) := \inf_{\vec{v}: || \vec{v} - \vec{1} || \le 1}  \left[ \ \frac{g(\vec{1}) - g(\vec{v})}{ ||\vec{v} - \vec{1}||^2} \ \right], \eqno(5.1)
$$
then $  \  \lambda(g) \in (0, \infty),  \  $ and as ordinary $ \ ||\vec{v}||^2 = ||v||^2  =  \sum_{j=2}^{n} v^2(j),  $ so that

$$
g(\vec{1}) - g(\vec{v}) \ge \lambda(g)  \cdot ||\vec{v} - \vec{1}||^2,  \  || \vec{v} - \vec{1} || \le 1. \eqno(5.2)
$$

 \ Further, let $  \ \epsilon \in (0, \lambda(g)) $  and denote also

$$
H_n(\lambda) := \sup_{ ||\vec{v} - \vec{1}||^2 \le \epsilon/\lambda } \left[  \prod_j |v(j)| \cdot h(\vec{v})  \right]; \eqno(5.3)
$$
then

$$
Q_n(\epsilon) \le H_n(\lambda) \cdot  \int_{ ||\vec{v} - \vec{1}||^2 \le \epsilon/\lambda  } \ d \vec{v} =
$$

$$
 H_n(\lambda) \cdot \frac{\pi^{(n-1)/2}}{\Gamma((n+1)/2)}  \cdot \left(  \frac{\epsilon}{\lambda(g)} \right)^{(n-1)/2}.  \eqno(5.4)
$$

 \vspace{4mm}

 \ The lower bound for this probability may be obtained quite analogously. Denote

$$
\mu = \mu(g) := \sup_{\vec{v}: || \vec{v} - \vec{1} || \le 1}  \left[ \ \frac{g(\vec{1}) - g(\vec{v})}{ ||\vec{v} - \vec{1}||^2} \ \right], \eqno(5.5)
$$
then $  \  \mu(g) \in (0, \infty),  \  $ and

$$
g(\vec{1}) - g(\vec{v}) \le \mu(g)  \cdot ||\vec{v} - \vec{1}||^2,  \  || \vec{v} - \vec{1} || \le 1. \eqno(5.6)
$$

 \  Let $  \ \epsilon \in (0, \mu(g)) $  and set also

$$
G_n(\mu) := \inf_{ ||\vec{v} - \vec{1}||^2 \le \epsilon/\mu } \left[  \prod_j |v(j)| \cdot h(\vec{v})  \right]; \eqno(5.7)
$$
then

$$
Q_n(\epsilon) \ge G_n(\mu) \cdot  \int_{ ||\vec{v} - \vec{1}||^2 \le \epsilon/\mu  } \ d \vec{v} =
$$

$$
 G_n(\mu) \cdot \frac{\pi^{(n-1)/2}}{\Gamma((n+1)/2)}  \cdot \left(  \frac{\epsilon}{\mu(g)} \right)^{(n-1)/2}.  \eqno(5.8)
$$

\vspace{4mm}

\section{ Examples. Concluding remarks.}

\vspace{4mm}

\ {\bf A.  An example.}  It is easily to verify that all the conditions of our theorem 2.1  are satisfied  for example for
arbitrary non-degenerate Normal (Gaussian) multivariate distribution, as well as  in the case when the random
variables $  \ \xi(i) \ $  are independent and have non-zero continuous density of distribution. \par

\vspace{4mm}

 \ {\bf B.  Possible generalizations.}  The offered here method may be easily generalized by our opinion on the asymptotic
computation for the distribution of a form

$$
{\bf P} \left( U(\vec{\xi})  > \max_{\vec{x}} U(\vec{x})  - \epsilon \right), \ \epsilon \to 0+,
$$
as well as for computation of integrals of a form

$$
I(\epsilon) = \int_{ y:  U(\vec{y}) > \max_{\vec{x}} U(\vec{x}) - \epsilon }   Z(y) \ \mu(dy)
$$
etc. \par

\vspace{4mm}

 \ {\bf C. Counterexamples.}  Let us prove that the condition about the existence of density function is essential for
our conclusions.  \par

\vspace{3mm}

 \ {\bf Example 6.1.}  Let the r.v.  $ \ \xi(j), \ j = 2, 3, \ldots, n-1  \ $ be arbitrary non-degenerate,  say,
 independent and  have the standard Gaussian distribution, and put $  \ {\bf P} (\xi(1) = 0) = 1. \ $ Then
both the r.v. $  \ \sum_j \xi(j), \ \sum_j \xi^2(j) $  have infinite differentiable bounded densities, but

$$
\sup_{\xi(i)} T(n) = \sqrt{n-1},
$$
therefore for sufficiently small positive values $  \ \epsilon $

$$
Q_n(\epsilon) = {\bf P} (T(n) > \sqrt{n} - \epsilon) = 0,
$$
in contradiction to the propositions of theorem 2.1 and 3.1. \par

\vspace{4mm}

 \ {\bf Example 6.2.}  Let now  the r.v.  $ \ \xi(i), \ i = 1, 2, 3, \ldots, n  \ $ be the Rademacher sequence, i.e. the sequence of
independent r.v. with distribution

$$
{\bf P}(\xi(i) = 1) = {\bf P}(\xi(i) = - 1) = 1/2.
$$

 \ Then

$$
{\bf P} (T(n) = \sqrt{n})  = 2^{-n} = {\bf P} (T(n) > \sqrt{n} - \epsilon),  \ 0 < \epsilon < (2 \sqrt{n})^{-1},
$$
in contradiction to the propositions of theorem 2.1 and 3.1. \par

 \vspace{8mm}

 {\bf References.}

 \vspace{4mm}

{\bf 1. Caballero M.E., Fernandes B. and Nualart D. }  {\it Estimation of densities
and applications.}  J. of Theoretical Probability, 1998, 27, 537-564.

\vspace{4mm}

{\bf 2. Xiequan Fan.} {\it Self-normalized deviations with applications to t-statistics.} \\
arXiv 1611.08436 [math. Pr] 25 Nov. 2016.

\vspace{4mm}

{\bf 3. Gine E., Goetze F. and Mason D. } (1997). {\it When is the Student t-
statistics asymptotically standard normal?}  Ann. Probab., 25, (1997), 1514-1531.

\vspace{4mm}

{\bf 4. I.Grama, E.Haeusler. }  {\it Large deviations for martingales. } Stoch. Pr.
Appl., 85, (2000), 279-293.

\vspace{4mm}

{\bf 5. B.Y.Jing, H.Y.Liang, W.Zhou.}  {\it Self-normalized moderate deviations
for independent random variables.}  Sci. China Math., 55 (11), 2012, 2297-2315.

\vspace{4mm}

{\bf 6. Ostrovsky E.I.} (1999). {\it Exponential estimations for Random Fields and
its applications,} (in Russian). Moscow-Obninsk, OINPE.

\vspace{4mm}

 {\bf 7. De La Pena V.H. }  {\it A general class of exponential inequalities for martingales and ratios.} 1999,
Ann. Probab., 36, 1902-1938.

\vspace{4mm}

{\bf 8. De La Pena V.H., M.J.Klass, T.L.Lai.} {\it Self-normalized Processes:
exponential inequalities, moment Bounds and iterative logarithm law.} 2004 ,
Ann. Probab., 27, 537-564.

\vspace{4mm}

{\bf 9. Q.M.Shao.} {\it Self-normalized large deviations.} Ann. Probab., 25, (1997),
285 - 328.

\vspace{4mm}

{\bf 10. Q.M.Shao. } {\it Self-normalized limit theorems: A survey.} Probability
Surveys, V.10 (2013), 69-93.

\vspace{4mm}

{\bf 11. Q.Y.Wang, B.Y.Jing.} {\it An exponential non-uniform Berry-Essen for
self-normalized sums.}  Ann. Probab., 27, (4), (1999), 2068-2088.

\end{document}